\documentclass[10pt]{elsart}

\journal{}

\usepackage[a4paper, left=0.8in,right=0.8in,top=0.8in,bottom=0.8in]{geometry}
\usepackage{amsmath}
\usepackage{amssymb}
\usepackage{leftidx}
\usepackage{url}

\begin{document}

\newtheorem {Theorem}{Theorem}
\newtheorem {Corollary}{Corollary}

\newcommand {\dd}{\textup{d}}
\newcommand{\RR}{\mathbb{R}}
\newcommand {\sts}[2]{\genfrac\{\}{0pt}{}{#1}{#2}}

\newcommand*\pFqskip{8mu}
\catcode`,\active
\newcommand*\W{\begingroup
        \catcode`\,\active
        \def ,{\mskip\pFqskip\relax}%
        \dopFq
}
\catcode`\,12
\def\dopFq#1#2#3{%
        W\left(\genfrac..{0pt}{}{#1}{#2};#3\right)%
        \endgroup
}

\begin{frontmatter}

 \title{Some physical applications of generalized Lambert functions}
 \author[mezo]{Istv\'an Mez\H{o}}
 \address[mezo]{Department of Mathematics,\\Nanjing University of Information Science and Technology,\\No.219 Ningliu Rd, Nanjing, Jiangsu, P. R. China}
 \ead{istvanmezo81@gmail.com}
 \thanks[mezo]{The research of Istv\'an Mez\H{o} was supported by the Scientific Research Foundation of Nanjing University of Information Science \& Technology, and The Startup Foundation for Introducing Talent of NUIST. Project no.: S8113062001}
 \author[keady]{Grant Keady}
 \address[keady]{Department of Mathematics, Curtin University,
 Perth, Australia}
 \ead{grant.keady@curtin.edu.au}
 
\begin{abstract}In this paper we review the physical applications of the generalized Lambert function recently defined by the first author. Among these applications we mention the eigenstate anomaly of the $H_2^+$ ion, the two dimensional two-body problem in general relativity,  the stability analysis of delay differential equations and water-wave applications. We also point out that the inverse Langevin function is nothing else but a specially parametrized generalized Lambert function. 
\end{abstract}

\end{frontmatter}

\section{Introduction}

\subsection{The definition of the Lambert $W$ function}
\label{subsec:defW}

The solutions of the transcendental equation
\[ xe^x=a \]
were studied by Euler and by Lambert : see~\cite{W}. The inverse of the function on the left hand side is called the Lambert function and is denoted by $W$. Hence the solution is given by $W(a)$. If $-\frac1e<a<0$, there are two real solutions, and thus two real branches of $W$ \cite{Veberic}. If we allow complex values of $a$, we get many solutions, and $W$ has infinitely many complex branches \cite{W,CJ,CJK,JHC}. 
These questions are comprehensively discussed by Corless et al.~\cite{W}. The question why 
we use the letter ``W'' is discussed by Hayes \cite{Hayes}; 
see also \cite{CJK}.

For the remainder of this paper, all parameters will be assumed to be real, and our
concern will be real-valued functions of real variables.

The Lambert function appears in many physical and mathematical problems. 
Some recently discussed physical models, however, lead to equations giving $\exp(x)$ as a given rational function of $x$:  
the preceding equation is but a simple special case. 
In the following section we briefly describe some significant applications of $W$.
Then we turn to the generalized Lambert function and some applications.

\subsection{Applications of $W$}

The survey paper of Corless et al. \cite{W} describes a large number of applications of $W$. This function appears in the combinatorial enumeration of trees, in the jet fuel problem, in enzyme kinetics, or in the solution of delay differential equations, just to mention a few areas. Some specific applications of $W$ in the investigation of solar wind comes from Cranmer \cite{Cranmer}, some applications in electromagnetic behavior of materials are given by Houari \cite{Houari}, other appearances in electromagnetics are investigated by Jenn \cite{Jenn}. It is known that Wien's displacement constant can be expressed by $W$, and in the discussion of capacitor fields $W$ also appears \cite{VJC}. The Lambert function has applications in quantum statistics, too \cite{VGJB}. A simple mathematical application connects $W$ to the distribution of primes via the Prime Number Theorem \cite{Visser}.

\section{The generalized Lambert function}

\subsection{Definition}
\label{subsec:defGenW}

Now we recall the recent generalization of the Lambert function and then we show its applicability in different physical problems where the classical $W$ function is insufficient.

In 2006 T. C. Scott \cite{SMM} and his co-workers defined the generalized $W$ function as the solution(s) of the equation
\[e^x\frac{(x-t_1)(x-t_2)\cdots(x-t_n)}{(x-s_1)(x-s_2)\cdots(x-s_m)}=a.\]
In  \cite{MezoBaricz}, where this generalization was studied, this generalized Lambert function is denoted by
\[\W{t_1,t_2,\dots,t_n}{s_1,s_2,\dots,s_m}{a}.\]

So, in particular,
\[\W{}{}{a}=\log(a),\quad\W{0}{}{a}=W(a),\quad\W{}{0}{a}=-W\left(-\frac{1}{a}\right)\]
\[ \W{t}{}{a}=t+W(a e^{-t}),\quad\W{}{s}{a}=s-W\left(-\frac{e^s}{a}\right). \]

We now describe various applications requiring the generalized Lambert function.

\subsection{Applications}
\label{subsec:Applications}

\subsubsection{A problem in molecular physics}\label{Problem1}

Corless and his co-workers \cite{W} mentioned that there was an anomaly in molecular physics. When physicists tried to calculate the eigenstates of the hydrogen molecular ion $(H_2^+)$, the results were not matching with the predictions \cite{Frost}. The problem was that -- being unable to solve their equations -- the physicists used numerical approximations which were inadequate. This problem originally emerged in 1956 and was solved by Scott and his coworkers \cite{SBDM} just in 1993. In the analytic solution the Lambert function appears which helped to take exponentially subdominant terms into account in the solution which could explain the anomaly. Briefly, using the double well Dirac delta function model, the wave equation to be solved was the following:
\[-\frac12\frac{d^2\psi}{dx^2}-q[\delta(x)+\lambda\delta(x-R)]\psi=E(\lambda)\psi.\]
The solution then becomes
\[\psi=Ae^{-d|x|}+Be^{-d|x-R|}.\]
Looking for the possible values of $d$, we make $\psi$ be continuous on each well, and this requirement gives that $d_\pm$ must satisfy the transcendental equation
\[d_{\pm}=q(1\pm e^{-d_\pm R}).\]
This can be converted to an equation solvable directly by using $W$, and we have that
\[d_\pm=q+\frac{W(\pm qRe^{-qR})}{R}.\]
Finally, the eigenenergies of the system are
\[E_{\pm}=-\frac{d_\pm^2}2.\]
This already matches to the predictions.

If one were to solve a similar eigenvalue problem for ions with unequal charges, $W$ is no longer sufficient. 
Instead, one must solve a transcendental equation of the form~\cite{SBDM}
\begin{equation}
e^{-cx}=a_0(x-t_1)(x-t_2).\label{geneq1}
\end{equation}
The solutions, as it is easy to see, can be given by the generalized Lambert function:
\[c\W{ct_1,ct_2}{}{\frac{c^2}{a_0}}.\]

\subsubsection{A problem in general relativity}\label{Problem2}

``One of the oldest and most notoriously vexing problems in gravitational theory is that of determining the (self-consistent) motion of $N$ bodies and the resultant metric they collectively
produce under their mutual gravitational influence\dots'' \cite{MO}. 
Mann and Ohta were considering the equations of motion of two bodies in one spatial dimension (to make the field equations tractable and eliminate gravitational waves, this dimensional restriction is very useful). The solution contains the $W$ function for the case when the two bodies have equal masses. When the masses are different, one needs to solve a more general equation of the form \cite{MO}
\begin{equation}
e^{-cx}=a_0(x-t_1)(x-t_2).\label{geneq1a}
\end{equation}

Refering back to the previous problem on the eigenstates of the hydrogen ion, we see that we have the identical mathematical problem.

\subsubsection{Bose-Fermi mixtures}

Scott mentioned \cite{SBDM} that in some recent physical investigations on Bose-Fermi mixtures the transcendental equation
\begin{equation}
e^{-cx}=a_0\frac{(x-t_1)(x-t_2)\cdots(x-t_n)}{(x-s_1)(x-s_2)\cdots(x-s_m)}\label{defmw}
\end{equation}
appears. We see that
\[c\W{ct_1,ct_2,\dots,ct_n}{cs_1,cs_2,\dots,cs_m}{\frac{c^{n-m}}{a_0}}\]
gives the solution(s) of equation \eqref{defmw}.

Independently of our study on the generalized $W$ function \cite{MezoBaricz}, Scott et al. \cite{SFG,SFGZh} studied the number of solutions of equations of the form \eqref{defmw}.

\subsubsection{Delay differential equations}

The generalized Lambert function helps in the solution of 
linear constant-coefficient delay differential equations
(and stability analysis for delay differential equations).

Consider linear constant-coefficient delay differential equations of the form
\[\sum_{k=0}^na_ku^{(k)}(t)=\sum_{k=0}^mb_ku^{(k)}(t-\tau).\]
Seeking solutions of the form $u(t)=\exp(\lambda t)$
leads to the characteristic equation
\[\sum_{k=0}^n a_k\lambda^k=e^{-\lambda\tau}\sum_{k=0}^m b_k\lambda^k.\]
Next, suppose the polynomials are factorized to give
\[ a_n \prod_{k=1}^n (\lambda - t_k) = b_m e^{-\lambda\tau}\prod_{k=1}^m (\lambda - s_k) .\]
This is an equation solvable by the generalized Lambert $W$ function.

To show an application, we cite S. A. Campbell: ``Second-order delay differential equations arise in a variety of mechanical, or neuro-mechanical systems in which inertia plays an important role\dots Many of these systems are regulated by feedback which depends on the state and/or the derivative of the state. In this
case the model equations take the form
\[\ddot u(t)+b\Dot u(t)+au(t)=f(u(t-\tau),\Dot u(t-\tau)),\]
where $a$, $b$ are positive constants representing physical attributes of the system, $\tau$ is the time delay, $u$, $u(t-\tau)$ are the values of the regulated variable evaluated at, respectively, times $t$ and $t-\tau$ and the function, $f(x, y)$, describes the feedback.'' \cite{Campbell}.
Linearizing about an equilibrium solution gives a linear constant-coefficient delay differential equations
(with $n=2$, $a_2=1$ and $m=1$ in the above) whose characteristic equation is
\begin{equation}
 (\lambda-t_1) (\lambda-t_2) = b_1 e^{-\lambda\tau} (\lambda - s_1) .
 \label{eq:dde2}
 \end{equation}


Another application arises in the paper of Chambers \cite{Chambers}. He considered some differential-difference equations with characteristic equations solvable by the generalized Lambert function with one upper and one lower parameter \cite[eq. (2.1b)]{Chambers}.

\subsubsection{Inverting dispersion equations for water waves}\label{subsubsec_wwave}

The dispersion relation for periodic water waves propagating into still water is 
\[
\omega^2 = g k\tanh\left(k h\right) ,
\]
where $\omega$ denotes the wave's frequency, $k$ its wave number,
$g$ the acceleration due to gravity and $h$ the mean water depth.
There are circumstances, e.g. in connection with computations using the
`Mild Slope Equation', where $g$, $h$ and $\omega$ are given and one needs
to solve for $k$. 
The numerical approximations for doing this are satisfactory, but it seems
worth noting that exact inversion is possible here, and in more general settings.

With $x=k h$ and $y={\omega^2 h}/{g}$ 
the preceding dispersion relation is
$
y = x \tanh(x)
$.
As $x \tanh(x) $ is increasing in $x$ for $x>0$, for each $y>0$ there is a unique $x>0$ solving this equation. 
Some simple bounds on the positive solution $x(y)$ are immediate.
As $x \tanh(x)<{\rm min}(x,x^2)$, $x>{\rm max}(y,\sqrt{y})$.
The equation rewrites to
\begin{equation}
e^{2 x} =\frac{x+y}{x-y} .
 \label{eq:MB6wwaves}
\end{equation}
and its solution is
\[
x = \frac{1}{2} 
W\left(\begin{array}{c} 2y\\-2y \end{array}; 1\right) .
\]

Waves in flows with two layers of different densities have been studied and
in~\cite[p. 421]{MT53} 
$$\omega^2 = g k\frac{\rho_1-\rho_2}{\rho_1\coth(k h)+\rho_2} . $$
Using the same definitions of $x$ and $y$ as in the preceding paragraph the equation above becomes
$$y= x\frac{\rho_1-\rho_2}{\rho_1\coth(x)+\rho_2} ,$$
which rearranges to
$$e^{2 x}= \frac{x+y}{x-\frac{\rho_1+\rho_2}{\rho_1-\rho_2} \, y} . $$
This becomes equation~(\ref{eq:MB6wwaves}) when $\rho_2=0$.
Once again one can find an expression for $x(y)$ in terms of generalized Lambert functions.



\subsubsection{The inverse Langevin function}

The Langevin function $L$ plays an important role in the study of paramagnetic materials and of polymers like rubber. More precisely, $L$ describes the magnetization of a paramagnet under the presence of outer classical (i.e., non quantical) magnetic field \cite{Kittel}. The expression of $L$ is
\[L(x)=\coth(x)-\frac1x.\]
The inverse of $L$ is widely studied, and good methods of approximating it are available. 
A comprehensive article on the approximations of $L^{-1}$ is \cite{Jedynak}.

We point out here that the inverse Langevin function $L^{-1}$ is also a special case of the generalized Lambert function. 
This can easily be seen by the fact that the transcendental equation $L(x)=a$ can be rewritten as
\[L(x)=\frac{e^{2x}+1}{e^{2x}-1}-\frac1x=a.\]
After some rearrangement and the variable substitution $x\to-\frac12x$ we get that this equation is equivalent to the equation
\begin{equation}
e^{-x} = \frac{a+1}{a-1} \frac{x-\frac{2}{a+1}}{x-\frac{2}{a-1}}.
\label{eq:invLang}
\end{equation}
Its solution is then given by
\[\W{\frac{2}{a+1}}{\frac{2}{a-1}}{\frac{a-1}{a+1}}.\]
Hence
\[L^{-1}(a)=-2\W{\frac{2}{a+1}}{\frac{2}{a-1}}{\frac{a-1}{a+1}}.\]



{\mbox{\rm or\ equivalently\ }} 

In what follows we collect some information from \cite{MezoBaricz} on some particular cases of the generalized Lambert function.

\section{The case of one upper and one lower parameter}
\label{sec:oneuponelow}

\begin{thm}\label{TaylorWts}The Taylor series of $\W{t}{s}{a}$ around $a=0$ is
\[\W{t}{s}{a}=t-T\sum_{n=1}^\infty\frac{L_n'(nT)}{n}e^{-nt}a^n,\]
where $T=t-s\neq 0$, and $L_n'$ is the derivative of the $n$th Laguerre polynomial \cite{Bell,Szego}.

Moreover, when $t<s$ the radius of convergence of the power series in Theorem \ref{TaylorWts} is the following
$$r=e^{t}\lim_{n\to\infty}\left|\frac {L_{n-1}^{(1)}(nT)} {L_n^{(1)}((n+1)T)}\right|=e^{\frac{t+s}{2}-2\sqrt{s-t}}.$$
\end{thm}

\subsection{The $r$-Lambert function}
\label{subsec:rLambert}

Mez\H{o} and Baricz \cite{MezoBaricz} studied the solution of the equation
\[xe^x+rx=n.\]
This can be rewritten as
\begin{equation}
-r = e^x \frac{x}{x -\frac{n}{r}},
\label{rLamdefeq}
\end{equation}
from where it is apparent that the solution is
\begin{equation} \W{0}{n/r}{-r}.
\label{WrWconn}
\end{equation}
This particular (generally multivalued) function is called $r$-Lambert function by the authors in 
\cite{MezoBaricz} and was denoted by $W_r(n)$.
It can be shown that
\[\lim_{r\to0}\W{0}{n/r}{-r}=W_0(n)=W(n),\]
the classical Lambert function. 
We also remark that the defining equation \eqref{rLamdefeq} of $W_r(n)$ does not contain any singularity around $r=0$.


Depending on the parameter $r$, the $r$-Lambert function has one, two or three real branches and so the above equations can have one, two or three solutions (we restrict our investigation to the real line). 

\begin{thm}\label{TaylorrLambert}Let $r\neq-1$ be a fixed real number. Moreover, we define the polynomial $M_k^{(n)}(y)$ as
\begin{equation}
M_k^{(n)}(y)=\sum_{i=1}^kn^{\overline i}\sts{k}{i}(-y)^i.\label{Mkdef}
\end{equation}
Then in a neighbourhood of $x=0$,
\begin{equation}
W_r(x)=\frac{x}{r+1}+\sum_{n=2}^\infty M_{n-1}^{(n)}\left(\frac{1}{r+1}\right)\frac{x^n}{(r+1)^nn!}.\label{TaylorWr}
\end{equation}
\end{thm}

It is known \cite{W} that for the classical Lambert $W$ function
\[W(x)\sim\log(x)+\log\left(\frac{1}{\log(x)}\right)\]
as $x\to\infty$. The asymptotics of the $r$-Lambert function is studied in \cite{MezoBaricz}: for any $r\in\RR$
\[W_r(x)\sim\log(x)+\log\left(\frac{1}{\log(x)}-\frac{r}{x}\right)\]
as $x\to\infty$. Moreover,
\[W_r(x)\sim\frac1rx\]
as $x\to-\infty$.

\section{The case of two upper parameters}

The Taylor series of $\W{t}{s}{a}$ (see Theorem \ref{TaylorWts}) contained easily identificable polynomials: the derivative of the Laguerre polynomials. In the case of $\W{t_1,t_2}{}{a}$ we have that the Bessel polynomials \cite{KF}
\[B_n(z)=\sum_{k=0}^n\frac{(n+k)!}{k!(n-k)!}\left(\frac{z}{2}\right)^k\]
appear (sometimes $y_n(z)$ notation is used). The following is a result of Mugnaini \cite{Mugnaini}.

\begin{thm}\label{TaylorWtt}The Taylor series of $\W{t_1,t_2}{}{a}$ around $a=0$ is
\[\W{t_1,t_2}{}{a}=t_1-\sum_{n=1}^\infty\frac{1}{n n!}\left(\frac{ane^{-t_1}}{T}\right)^nB_{n-1}\left(\frac{-2}{nT}\right),\]
where $T=t_2-t_1$.
\end{thm}

\section*{Closing remarks}

The first author wrote C code to calculate the real $r$-Lambert function on all the branches. This can be downloaded from \url{https://sites.google.com/site/istvanmezo81/}.

\end{document}